\documentclass{article}

\usepackage{arxiv}

\usepackage[english]{babel}
\usepackage{biblatex}
\usepackage{amsmath}
\usepackage[pdftex]{hyperref}
\usepackage[utf8]{inputenc} % allow utf-8 input
\usepackage[T1]{fontenc}    % use 8-bit T1 fonts
\usepackage{hyperref}       % hyperlinks
\usepackage{url}            % simple URL typesetting
\usepackage{booktabs}       % professional-quality tables
\usepackage{amsfonts}       % blackboard math symbols
\usepackage{nicefrac}       % compact symbols for 1/2, etc.
\usepackage{microtype}      % microtypography
\usepackage{lipsum}
\usepackage{graphicx}
\graphicspath{ {./images/} }

\title{General formula to solve quintic equation}

\newmuskip\pFqmuskip

\newcommand*\pFq[6][8]{%
  \begingroup % only local assignments
  \pFqmuskip=#1mu\relax
  \mathchardef\normalcomma=\mathcode`,
  \mathcode`\,=\string"8000
  \begingroup\lccode`\~=`\,
  \lowercase{\endgroup\let~}\pFqcomma
  {}_{#2}F_{#3}{\left[\genfrac..{0pt}{}{#4}{#5};#6\right]}%
  \endgroup
}
\newcommand{\pFqcomma}{{\normalcomma}\mskip\pFqmuskip}

\author{
 Rodrigo José Martinelli Biglia Andrade \\
  Administration Department\\
  Methodist University of Piracicaba\\
  Piracicaba, SP, Brazil \\
  \texttt{contato@pontoabc.com} \\
}

\begin{document}
\maketitle
\begin{abstract}
There are many papers about the proof of Abel-Ruffini theorem, but, in this paper, we will show a way to solve any quintic equation by radicals which will be a outstandig result because of work the mathematicians Abel and Ruffini who showed that there is no formula to quintic equation. We need to mention the mathematician Évariste Galois that also showed the same conclusion of Abel and Ruffini using permutations,group theory etc. However, in this article we will present a method to split any quintic equation into two equation (one of degree 3 and another of degree 2) where both equations can be solved by radicais using the quadratic and the cubic formula.
\end{abstract}

% keywords can be removed

\section{Introduction}
This paper contributes to a new way to solve quintic equation using a polynomial (Martinelli's polynomial) that can be proved using some simple steps. Since the 1500's when the cubic and quartic formula was discovered the world waited centuries for the next step: to find a method or general formula to solve the quintic equation. Although Abel and Ruffini showed the impossibility of a closed formula to solve general quintic equation the search for a formula to solve quintic equation ends. But the mathematician Joseph-Louis Lagrange proved the inversion theorem that could solve specific cases of the quintic equation and higher degree. However using hypergeometric function will be possible to solve any quintic equation if reducing the general quintic equation to the Bring-Jerrard form (but it is a very complex process to transform the general quintic equation to that form). At the end of this paper there is an example that checks if the method that will be presented works. Finally, to end the abstract we will not proof or talk about the Lagrange inversion theorem and Bring-Jerrard form to get more understanding this paper.
\section{About Galois theory and Abel-Ruffini theorem}
\label{sec:headings}
We will see some concepts about Galois' theory and Abel-Ruffini's theorem.
The theorem, in a succinct way, consists of the proof that there is no "closed" formula for all equations of degree greater or equal to 5. That is, there is no way to algebraically arrive at a formula for those equations greater or equal to 5. As Zoladek (2000) reported, "A general algebraic equation of degree $\geq$ 5 cannot be solved in radicals. This means that there does not exist any formula
which would express the roots of such equation as functions of the coefficients by means of the algebraic operations and roots of natural degrees."  (p. 254).

Thus, Galois theory comes to confirm, through the study of the relationship of the roots that it is impossible to solve, algebraically, equations of degree equal to or greater than 5. In its most basic form, the theorem states that given a
E / F field extension which is finite and Galois, there is a one-to-one correspondence between its intermediate fields and subgroups of its Galois group. (Intermediate fields are K fields that satisfy $F \subseteq K \subseteq E$; they are also called sub-extensions of E / F.).

So, according to Abel-Ruffini's theorem and Galois theory, there are cases where equations of degree 5 or higher are solvable by radicals or algebraically.

\section{Proof of Martinelli's polynomial}

To demonstrate the Martinelli's Polynomial it is necessary to consider the roots of any fifth degree equation.
Each root must be associated with a single other root. We have then that the combination of the roots of a fifth degree equation is ten (due that the polynomial is tenth degree). Let us then consider the following equations:

\begin{equation}
(x^2-kx+n)(x^3+kx^2+kx+m)=0.
\end{equation}
\begin{equation}
(x^2-kx+n)(x^3+kx^2+lx+m)=0.
\end{equation}

Equation (1) have a relation with equation (2). If we expand both equations (1) and (2), each corresponding term can be factorized. Let put those terms into a table:

\def\arraystretch{2.0}%
\begin{table}[htb]
\centering
\caption{Equation (1) expanded}
\begin{tabular}{|c|c|c|c|}
\hline
$x^3$ & $x^2$ & $x$ & Term Independent \\
\hline
$k+n-k^2=C_2$ & $-k^2+kn+m=D_2$ & $kn-km=E_2$ & $ mn=F$ \\
\hline
\end{tabular}
\end{table}
\def\arraystretch{2.0}%
\begin{table}[htb]
\centering
\caption{Equation (2) expanded}
\begin{tabular}{|c|c|c|c|}
\hline
$x^3$ & $x^2$ & $x$ & Term Independent \\
\hline
$l+n-k^2=C$ & $nk-lk+m=D$ & $nl-km=E$ & $ mn=F$ \\
\hline
\end{tabular}
\end{table}

Andrade (2019) consider the follow equation to better understand the proof of Martinelli's Polynomial:

\begin{equation}
x^2+3x+2=0
\end{equation}

If we plug in x the sum of the roots of equation (3) we have 2=0, that would be absurd, but with that idea, we can create a second degree equation using the terms $C, C_2 , D, D_2, E$ e $E_2$ as shown in Tables 1 and 2 above. Matching the equation that will be created with $(-k^2+n+k-C)n$ we will have, on the right side of the equation, the same thing $E_2- E$, because $E_2$ is the same as $kn-km$ and $-E$ is the same as $-nl+km$. The sum of $E_2-E$ will result $kn-nl$ then $l=C-n+k^2$. So we have $kn -n(C-n+k^2)$ that is the same as $(-k^2+n+k-C)n$.

To create a second degree equation where one of the solutions of the equation will be the sum of two solutions of a fifth degree equation we must follow this logic: We have the general form of the second degree equation that is
$ax^2 + bx + c = 0$, the coefficients of the second degree equation that will be created will be $a = C_2- C$, $b = D_2- D$, $c=E_2- E$ and on the right side of the equation we have to add $(-k^2+n+k-C)n$ because on the left side of the equation will remain $E_2- E$.

This idea will give us the possibility to know n from equation (2). The goal is to form a tenth degree equation with the unknown k. So $a = C_2 - C$ and $b=D_2-D$ according to tables 1 and 2, it follows that $a=k+n-k^2-C$, $b=-k^2+kn+m-D$, $c=nk -km-E$ and the right side of the equation will be $(-k^2+n+k-C)n$. Thus, making the substitutions in $ak^2+bk+c=0$ we will have:

\begin{align*}
(-k^2+n+k-C)k^2 + (-k^2+nk+m-D)k + nk-km-E &= (-k^2+n+k-C)n.\\
(-k^4+nk^2+k^3-Ck^2)+(-k^3+nk^2+km- Dk)+ nk-km- E &= (-k^2+n+k-C)n. \\
-k^4+2nk^2-Ck^2-Dk+ nk- E &= -nk^2+n^2+nk-nC.\\
-k^4-Ck^2-Dk-E&=-3nk^2+n^2-nC.\\
\end{align*}

Putting n on the left side of the equation, we have:

\begin{equation}
n = \frac{k^4+Ck^2+Dk+E}{3k^2-n+C}.
\end{equation}

Since k is the sum of two roots of a fifth degree equation, C, D, E and F the coefficients of a fifth degree equation and n is the product of the roots of a second degree equation, so to arrive at Martinelli's polynomial, we have to replace the variable n with the use of an algebraic manipulation. This algebraic manipulation consists of taking the equality referring to the term D in table (2). Thus:

\begin{equation}
nk-(k^2-n+C)k+\frac{F}{n}=D.
\end{equation}

\begin{equation}
n^2=\frac{nk^3+Cnk-F+Dn}{2k}.
\end{equation}

Now just replace $ n^2$ in equation (4) and get n:

\begin{equation*}
n = \frac{k^4+Ck^2+Dk+E}{3k^2-n+C}.
\end{equation*}
\begin{equation*}
-n^2+3nk^2+Cn=k^4+Ck^2+Dk+E.
\end{equation*}
\begin{equation*}
\frac{-nk^3-Cnk+F-Dn}{2k} + 3nk^2 + Cn = k^4+Ck^2+Dk+E.
\end{equation*}
\begin{equation}
n = \frac{2(k^5+Ck^3+Dk^2+Ek)-F}{5k^3+Ck-D}.
\end{equation}

If we have n, then we can create the Martinelli's polynomial that will be very important to solve any quintic equation
\begin{equation*}
\frac{2(k^5+Ck^3+Dk^2+Ek)-F}{5k^3+Ck-D} = \frac{k^4+Ck^2+Dk+E}{3k^2-\frac{2(k^5+Ck^3+Dk^2+Ek)-F}{5k^3+Ck-D}+C}.
\end{equation*}

Arranging the equation on both sides and equaling 0 we arrive at the Martinelli's polynomial:

\begin{eqnarray*}
{(2(k^5+Ck^3+Dk^2+Ek)-F)(13k^5+6Ck^3-5Dk^2+(-2E+C^2)k+F-DC)}
\end{eqnarray*}
\begin{equation}
-(k^4+Ck^2+Dk+E)(5k^3+Ck-D)^2 = 0.
\end{equation}
\section{Martinelli's Polynomial expanded}

\begin{equation*}
k^{10} + 3Ck^8 + Dk^7 + (3C^2-3E)k^6 + (2DC-11F)k^5 + (C^3-D^2-2CE)k^4 + (DC^2-4DE-4CF)k^3
\end{equation*}
\begin{equation}
 + (7DF-CD^2-4E^2+EC^2)k^2 + (4EF-FC^2-D^3)k-F^2+FDC-D^2E=0.
\end{equation}

\section{Solving a quintic equation as an example}

Let's begin with the quintic equation:

\begin{equation}
x^{5}+x+3=0
\end{equation}

Using the Martinelli's Polynomial (9) we get the follow polynomial:

\begin{equation}
k^{10}-3k^{6}-33k^{5}-4k^{2}+12k-9=0
\end{equation}

The roots of the equation (11) are the combination of the sums of the equation (10). So, wouldn't be necessary to know all roots of the equation (11). But, solving the equation (10) we will see that the sum of the complex roots give us one of the real root: 2.0837590792241645736... that is represented using hypergeometric function as follows:

\begin{equation}
k = \sqrt{2}\pFq{4}{3}{-\frac{1}{20},\frac{3}{20},\frac{7}{20},\frac{11}{20}}{\frac{1}{4},\frac{1}{2},\frac{3}{4}}{-\frac{253125}{256}}
-\frac{45\pFq{4}{3}{-\frac{9}{20},\frac{13}{20},\frac{17}{20},\frac{21}{20}}{\frac{3}{4},\frac{5}{4},\frac{3}{2}}{-\frac{253125}{256}}}{16\sqrt{2}}
+\frac{3}{2}\pFq{4}{3}{\frac{1}{5},\frac{2}{5},\frac{3}{5},\frac{4}{5}}{\frac{1}{2},\frac{3}{4},\frac{5}{4}}{-\frac{253125}{256}}
\end{equation}

\begin{equation}
\Bigg(x^2-kx+\frac{2(k^5+k)-3)}{5k^3}\Bigg)\Bigg(x^3+kx^2+\Bigg(k^2-\frac{2(k^5+k)-3)}{5k^3}\Bigg)x+\frac{15k^3}{2(k^5+k)-3}\Bigg)=0
\end{equation}

Now, plugging the k root, represented by hypergeometric function (12), we will have the equation (10), but, in this case, the equation (10) was factorized into two equations of degree 2 and another of degree 3. Using the quadratic and cubic formula we will able to solve in terms of radicals the equation (13) which it is the same as the equation (10).

\section{Remark}

As seen the Martinelli's polynomial provides the sum of two roots of any quintic equation and this is very important to split a quintic equation into two smaller degree equations. This method solves quintic equation using quadratic and cubic formulas but in the most cases it is necessary to use the hypergeometric function to represent the variable that splits the equation. However, a general method has been achieved to solve any quintic equation. But, it's no make sense, because we need to know one of the real root of the Martinelli's polynomial once that is the sum of quintic roots equation. We need to know before the roots of the quintic equation in the Bring-Jerrard form to split the fifth degree equation and solve in terms of radicals. So, Abel-Ruffini theorem and Galois theory didn't give us a clear answer about this way to solve quintic equation in terms of radicals and in the Bring-Jerrard form. Maybe it's a new discussion about equations theory field.

\end{document}